\documentclass[11pt,reqno]{amsart}
\usepackage{amsmath,amssymb,amsfonts,amsthm}
\usepackage[T1]{fontenc}
\usepackage[numbers,sort&compress]{natbib}
\usepackage{txfonts}
\usepackage[colorlinks=true,allcolors=blue]{hyperref}
\hypersetup{unicode=true,colorlinks=true,linkcolor=blue,citecolor=blue,urlcolor=blue,pdfauthor={Hirofumi Ota}}

\theoremstyle{plain}
\newtheorem{theorem}{Theorem}[section]
\newtheorem{proposition}[theorem]{Proposition}
\newtheorem{lemma}[theorem]{Lemma}

\theoremstyle{definition}

\newtheorem{example}[theorem]{Example}
\theoremstyle{remark}
\newtheorem{remark}[theorem]{Remark}

\newcommand{\R}{\mathbb{R}}
\newcommand{\E}{\mathbb{E}}
\newcommand{\Pbb}{\mathbb{P}}
\newcommand{\1}{\mathbf{1}}
\newcommand{\ph}{\varphi}
\newcommand{\ip}[2]{\langle #1,#2\rangle}
\newcommand{\norm}[1]{\lVert #1\rVert}
\newcommand{\abs}[1]{\lvert #1\rvert}
\newcommand{\eps}{\varepsilon}
\DeclareMathOperator{\diam}{diam}
\DeclareMathOperator{\dist}{dist}
\DeclareMathOperator{\conv}{conv}
\DeclareMathOperator*{\argmax}{argmax}
\DeclareMathOperator*{\argmin}{argmin}

\title[Near-Maximizers of Gaussian Processes]{Sharp Bounds for Near-Maximizers of Gaussian Processes}
\author{Hirofumi Ota}
\address{Komaba Institute for Science, Graduate School of Arts and Sciences, The University of Tokyo, 3-8-1 Komaba, Meguro-ku, Tokyo 153-8902, Japan}
\email{hirofumi-ota@g.ecc.u-tokyo.ac.jp}
\date{}
\hypersetup{pdftitle={Sharp Bounds for Near-Maximizers of Gaussian Processes}}

\begin{document}

\begin{abstract}
For a centered Gaussian process $G$ with continuous sample paths on a compact metric space $T$, a continuous drift $f:T\to\R$, a noise level $\sigma>0$, and $\eps\ge0$, let $A_\eps$ be the set of points at which $f+\sigma G$ is within $\eps$ of its maximum. Writing $\E\sup_TG$ for the expected supremum of $G$, we prove that the expected squared diameter of $A_\eps$ in the canonical pseudometric of $G$ is at most $8\eps \E\sup_TG/\sigma$ and that the expected maximal squared distance from a point of $A_\eps$ to an exact maximizer is at most $4\eps \E\sup_TG/\sigma$. The constants cannot be improved uniformly over index sets, processes, and drifts, and no logarithmic factor in $\eps$ appears, unlike in bounds integrated from tail estimates. At $\eps=0$ the bounds give the almost sure uniqueness of the maximizer in the canonical pseudometric without nondegeneracy assumptions. Given this uniqueness, both bounds follow from a bound with constant $2$ on the mean square of the largest increment of an independent copy of $G$ from an exact maximizer to a point of $A_\eps$. The proof rests on an identity for the mean squared support function of the $\eps$-subdifferentials of a smooth convex function at a standard Gaussian vector, obtained from two Gaussian integrations by parts.
\end{abstract}
\maketitle

\section{Introduction}
\label{sec:intro}

Consider a centered Gaussian process $G$ on a compact index set $T$, a deterministic function $f$ on $T$, and a noise level $\sigma>0$. Under mild conditions the maximizer of $f+\sigma G$ is almost surely unique~\cite{KimPollard1990,Arcones1992}. This paper gives a quantitative form of this fact. A point of $T$ is an $\eps$-near-maximizer if the value of $f+\sigma G$ at this point is within $\eps\ge0$ of the maximum. We bound the size of the set of $\eps$-near-maximizers in the canonical pseudometric of $G$. The bounds are uniform in $f$ and their constants are sharp.

\subsection{Main result}
\label{sec:mainresult}
Let $(T,\rho)$ be a nonempty compact metric space, and let $C(T)$ be the space of real-valued continuous functions on $T$. For any function $g:T\to\R$, write $\norm g_\infty=\sup_{t\in T}\abs{g(t)}\in[0,\infty]$; we equip $C(T)$ with this uniform norm and the associated Borel $\sigma$-field. On a probability space $(\Omega,\mathcal F,\Pbb)$, let $G$ be a Borel random element of $C(T)$ whose finite-dimensional distributions are centered Gaussian; then we have $\E\norm G_\infty<\infty$~\cite{AdlerTaylor2007}. For a deterministic drift $f\in C(T)$, a noise level $\sigma>0$, and $\eps\ge0$, define the process $X$, its set $A_\eps$ of $\eps$-near-maximizers, and the canonical pseudometric $d$ of $G$ by
\begin{align}
 X(t)&=f(t)+\sigma G(t),&
 A_\eps&=\{t\in T:X(t)\ge\max_T X-\eps\},\quad \eps\ge0,\label{eq:nearset}\\
 d(s,t)&=\bigl(\E(G(s)-G(t))^2\bigr)^{1/2}.&&\label{eq:canonical}
\end{align}
For nonempty $B\subset T$, write $\diam_d(B)=\sup_{s,t\in B}d(s,t)$ for its $d$-diameter. The sets $A_\eps$ are nonempty and compact, and $A_0$ is the set of exact maximizers of $X$.

\begin{theorem}\label{thm:process}
For every $\eps\ge0$, we have
\begin{align}
 \E\diam_{d}(A_\eps)^2
 &\le\min\Bigl\{\frac{8\eps}{\sigma}\E\sup_TG,\ \diam_d(T)^2\Bigr\},\label{eq:processbound}\\
 \E\max_{\substack{t\in A_\eps\\p\in A_0}}d(t,p)^2
 &\le\min\Bigl\{\frac{4\eps}{\sigma}\E\sup_TG,\ \diam_d(T)^2\Bigr\}.\label{eq:processradius}
\end{align}
Moreover, let $\widetilde G$ be an independent copy of $G$, and put
\[
 \Gamma_\eps=\max_{t\in A_\eps}\widetilde G(t)-\max_{p\in A_0}\widetilde G(p).
\]
Then, for every $\eps\ge0$, we have
\begin{equation}
 \E\Gamma_\eps^2\le\frac{2\eps}{\sigma}\E\sup_TG.
 \label{eq:processsharp}
\end{equation}
In particular, \eqref{eq:processbound} with $\eps=0$ gives $\diam_d(A_0)=0$ almost surely, so the maximizer is almost surely unique when $d$ is a metric.
\end{theorem}

The constants $8$, $4$, and $2$ are optimal (Proposition~\ref{prop:sharpness}). The bounds concern the whole set $A_\eps$, so they apply to any measurable selection from it, however it depends on $G$, provided $\eps$ is fixed in advance; the case of a random $\eps$ is discussed in Remark~\ref{rem:assumptions}.

The simplest example is Euclidean. Let $k\ge1$ be an integer, let $\ip\cdot\cdot$ and $\norm\cdot$ be the Euclidean inner product and norm on $\R^k$, let $I_k$ be the $k\times k$ identity matrix, and write $\diam$, without subscript, for the Euclidean diameter. If $T\subset\R^k$ is compact and $G(t)=\ip Zt$ with $Z\sim N(0,I_k)$, the canonical Gaussian process on $T$, then $d$ is the Euclidean distance and $\E\sup_TG$ is the Gaussian width of $T$. Here $T$ may be nonconvex or finite, and the numerical constants do not depend on $k$: the size of $T$ enters only through the Gaussian width, which for finite $T$ with $\abs T$ points is at most $\max_{t\in T}\norm t\sqrt{2\log\abs T}$.

\begin{example}\label{ex:brownian}
Let $G$ be a standard Brownian motion on $T=[0,1]$. Then we have $d(s,t)=\abs{s-t}^{1/2}$ and, by the reflection principle, $\E\sup_TG=\E\abs{G(1)}=\sqrt{2/\pi}$. For every $f\in C(T)$, $\sigma>0$, and $\eps\ge0$, \eqref{eq:processbound} therefore gives
\[
 \E\sup\{\abs{s-t}:s,t\in A_\eps\}\le\frac{8\eps}{\sigma}\sqrt{\frac2\pi}.
\]
The order $\eps$ is sharp already for $f=0$ and $\sigma=1$: since the increments of $G$ after time $3/4$ are independent of its earlier path, there is $c>0$ such that, for $0<\eps\le1$, we have
\[
 \Pbb\Bigl\{\min\Bigl(\max_{[0,1/4]}G,\ \max_{[3/4,1]}G\Bigr)\ge\max_TG-\eps\Bigr\}\ge c\eps.
\]
On this event $A_\eps$ meets both $[0,1/4]$ and $[3/4,1]$, so we have $\E\sup\{\abs{s-t}:s,t\in A_\eps\}\ge c\eps/2$.
\end{example}

\subsection{Comparison with previous bounds}
\label{sec:tails}

Block et al.~\cite[Theorem 5]{BlockEtAl2024}, sharpening~\cite{BlockEtAl2022}, prove the following tail bound in their analysis of Gaussian-perturbed empirical risk minimization: if $d$ is a metric and $G$ has unit marginal variances, then, in the present notation, it holds that
\begin{equation}
 \Pbb\left\{\max_{\substack{t\in A_\eps\\p\in A_0}}d(t,p)>r\right\}
 \le\frac{2\eps}{\sigma r^2}\E\sup_TG,
 \qquad r>0.
 \label{eq:priortail}
\end{equation}
Their result is stated for minimizers, in terms of a covariance normalized by the variance at the minimizer. For marginal variances in $[\kappa,1]$ with $0<\kappa<1$, the bound carries the factor $\kappa^{-2}$ and no longer translates into a bound on $d$; points of zero variance, such as $t=0$ for Brownian motion on $[0,1]$ or for $G(t)=\ip Zt$ when $0\in T$, are excluded.
Their proof treats one distance level at a time, which is why it yields a tail bound. Passing from \eqref{eq:priortail} to the second moment costs a logarithm, as in~\cite{BlockEtAl2022}: if $0<2\eps \E\sup_TG/\sigma\le\diam_d(T)^2$, integrating the smaller of one and the right-hand side of \eqref{eq:priortail} against $2r dr$ over $0<r\le\diam_d(T)$ gives
\[
 \E\max_{\substack{t\in A_\eps\\p\in A_0}}d(t,p)^2
 \le\frac{2\eps}{\sigma}\E\sup_TG\Bigl(1+\log\frac{\sigma\diam_d(T)^2}{2\eps \E\sup_TG}\Bigr).
\]
The logarithm counts distance scales: between $(2\eps \E\sup_TG/\sigma)^{1/2}$ and $\diam_d(T)$, each dyadic range of $r$ contributes the same amount. It cannot be removed using \eqref{eq:priortail} and boundedness alone, since a random variable with values in $[0,\diam_d(T)]$ whose tail probabilities are the smaller of one and the right-hand side of \eqref{eq:priortail} has this second moment.

Theorem~\ref{thm:process} has no logarithm because its proof integrates over $\eps$ instead. For a finite index set, written as $T\subset\R^k$ with $G(t)=\ip Zt$, two Gaussian integrations by parts express a regularized counterpart of $\E\Gamma_\eps^2$ as the integral, from $0$ to $\eps/\sigma$, of twice the mean of $\ip Zq$ at a random point $q$ that lies, up to a vanishing error, in the convex hull of $T$ (Lemma~\ref{lem:support}). Since $\ip Zq\le\sup_{t\in T}\ip Zt$ on this hull, the mean square grows by at most $2 \E\sup_TG$ per unit of $\eps/\sigma$, at whatever distance from the maximizer the near-maximizers lie. This gives \eqref{eq:processsharp} with its optimal constant and without a condition on the marginal variances. The tail bound itself is not improved: combined with Markov's inequality, \eqref{eq:processradius} gives \eqref{eq:priortail} with the constant $4$ in place of $2$.

For finite $T\subset\{0,1\}^k$ and $G(t)=\ip Zt$, the proof of~\cite[Lemma 6]{BuzagloHazan2026} contains a bound without a logarithm in $\eps$: when $T$ has $m$ points and $\eps\le\sigma$, the expected number of coordinates on which the points of $A_\eps$ do not all agree, which is at least $\E\diam(A_\eps)^2$, is at most of order $(\eps/\sigma)(k\log m)^{1/2}$. An identity in that proof implies $\limsup_{\eps\downarrow0}\eps^{-1} \E\diam(A_\eps)^2\le2 \E\sup_TG/\sigma$ for fixed $T$, a quarter of the constant in \eqref{eq:processbound}, but the argument conditions on one coordinate at a time and relies on the independence of the coordinates. The rotation argument of Chatterjee, Dembo, and Ding~\cite{ChatterjeeDemboDing2013} for level sets of Gaussian fields, adapted to $f+\sigma\gamma^{-1}(\gamma G+(1-\gamma^2)^{1/2}\widetilde G)$ with $\gamma^{-1}=1+\eps/(\sigma\E\sup_TG)$ for $\eps>0$ and $\E\sup_TG>0$, gives the first-moment bound $\E\Gamma_\eps\le(2\eps \E\sup_TG/\sigma)^{1/2}$ but does not seem to give \eqref{eq:processsharp}. Theorem~\ref{thm:process} also bounds the probability of a near-tie: if $S,S'\subset T$ are nonempty compact sets with $d(s,s')\ge r>0$ for all $s\in S$ and $s'\in S'$, then \eqref{eq:processradius} on the index set $S\cup S'$ and Markov's inequality give
\[
 \Pbb\Bigl\{\bigl|\max_SX-\max_{S'}X\bigr|\le\eps\Bigr\}\le\frac{4\eps}{\sigma r^2} \E\sup_{S\cup S'}G.
\]
If the variances of $G$ on $S$ are at least $r^2$, then adjoining to $S$ a point $s'$ with $G(s')=0$ and $f(s')=x$, and taking $S'=\{s'\}$, we obtain
\[
 \sup_{x\in\R}\Pbb\Bigl\{\bigl|\max_SX-x\bigr|\le\eps\Bigr\}\le\frac{4\eps}{\sigma r^2} \E\max\Bigl\{\sup_SG, 0\Bigr\},
\]
a bound on the L\'evy concentration function of $\max_SX$ that for equal variances has the form of the anti-concentration bound of~\cite[Theorem 3]{ChernozhukovChetverikovKato2015}. Finally, Theorem~\ref{thm:process} quantifies the almost sure uniqueness of the maximizer on which the argmax continuous mapping theorem relies~\cite{KimPollard1990,vanderVaartWellner1996}; Example~\ref{ex:chernoff} treats Brownian motion with parabolic drift on the line.

\subsection{Organization}
\label{sec:organization}
Section~\ref{sec:proof} proves an identity for the mean squared support function of the $\eps$-subdifferentials of a smooth convex function at a standard Gaussian vector (Lemma~\ref{lem:support}) and outlines how Theorem~\ref{thm:process} follows from it. Section~\ref{sec:sharpness} states the optimality of the constants and a drift-adaptive form of the bounds, and ends with a remark on the assumptions of Theorem~\ref{thm:process}. The appendices contain the remaining proofs.

\section{Proof of the main theorem}
\label{sec:proof}
For a finite index set, the Gaussian process is a linear functional of a standard Gaussian vector, and the vectors representing $\eps$-near-maximizers are $(\eps/\sigma)$-subgradients of a convex value function at that vector. Theorem~\ref{thm:process} thus reduces to bounding the size of $\eps$-subdifferentials of convex functions at a standard Gaussian vector. After a smooth strongly convex regularization, the identity \eqref{eq:supportidentity} expresses the mean squared support function of the $\eps$-subdifferential as an integral over $0<s\le\eps$ of the inner product of this Gaussian vector with the Gaussian average of the support points of the $s$-subdifferential, through two Gaussian integrations by parts. The Gaussian width bounds this inner product; this gives the constant $2$ in \eqref{eq:processsharp}, and Lemma~\ref{lem:geometry} then gives $8$ and $4$. This section proves the identity; the remaining steps are in Appendices~\ref{app:regularity}--\ref{sec:process}.

The ingredients are classical, although we have not found the identity in the literature. The support functions of $\eps$-subdifferentials are the $\eps$-directional derivatives studied in~\cite{Auslender1982,HiriartUrrutyLemarechal1993}. Gaussian integration by parts with respect to this Gaussian vector underlies the trace bound in~\cite{AbernethyEtAl2014}. The Gaussian average of the support points is the Steiner point of the $s$-subdifferential, in the Gaussian form of~\cite{Vitale1996Covariance}.

Throughout, $\ph(x)=(2\pi)^{-1/2}\mathrm e^{-x^2/2}$ denotes the standard normal density, $k\ge1$ denotes the Euclidean dimension, $a_+$ denotes $\max\{a,0\}$ for $a\in\R$, and $\1$ denotes an indicator.

\subsection{A support-function identity}
\label{sec:support}

\begin{lemma}\label{lem:geometry}
Let $B\subset\R^k$ be compact with $0\in B$, let $h_B(u)=\max_{x\in B}\ip ux$, $u\in\R^k$, be its support function, and let $U\sim N(0,I_k)$. Then we have
\begin{equation}
 \diam(B)^2\le4\E h_B(U)^2,\qquad
 \max_{x\in B}\norm x^2\le2\E h_B(U)^2.
 \label{eq:geometry}
\end{equation}
\end{lemma}
\begin{proof}
Choose $x,y\in B$ attaining the diameter. Since $\abs{\ip u{x-y}}\le h_B(u)+h_B(-u)$ for all $u\in\R^k$, we have
\[
 \norm{x-y}^2=\E\ip U{x-y}^{2}
 \le\E(h_B(U)+h_B(-U))^2
 \le2\E(h_B(U)^2+h_B(-U)^2)=4\E h_B(U)^2.
\]
For $x\in B$ of maximal norm, we have $(\ip Ux)_+\le h_B(U)$, because $h_B\ge0$ when $0\in B$. By the symmetry of the Gaussian law, we have $\norm x^2=2\E(\ip Ux)_+^2\le2\E h_B(U)^2$.
\end{proof}

For the next lemma, let $H:\R^k\to\R$ be twice continuously differentiable and $\mu$-strongly convex for some $\mu>0$, that is, $u^\top\nabla^2H(t)u\ge\mu\norm u^2$ for all $t,u\in\R^k$. A superscript $*$ denotes convex conjugation, defined by
\[
 H^*(z)=\sup_{t\in\R^k}\{\ip zt-H(t)\},
 \qquad H^{**}=(H^*)^*.
\]
For $z\in\R^k$ and $s\ge0$, define
\begin{equation}
 \begin{aligned}
 p(z)&=\argmin_{t\in\R^k}\{H(t)-\ip zt\},\\
 C_s(z)&=\{t\in\R^k:H(t)-\ip zt+H^*(z)\le s\}.
 \end{aligned}
 \label{eq:smoothdefs}
\end{equation}
The expression $H(t)+H^*(z)-\ip zt$ is the Fenchel--Young gap. Since $H^{**}=H$, the set $C_s(z)$ is the $s$-subdifferential of $H^*$ at $z$, that is, the set $\{t:H^*(z')\ge H^*(z)+\ip t{z'-z}-s\ \text{for all }z'\in\R^k\}$. Indeed, both conditions are equivalent to $\sup_{z'}\{\ip t{z'}-H^*(z')\}\le\ip tz-H^*(z)+s$. For $s\ge0$ and $z,u\in\R^k$, put
\begin{equation}
 \begin{aligned}
 q_s(z,u)&=
 \begin{cases}
  \argmax_{t\in C_s(z)}\ip ut,&s>0,\ u\ne0,\\
  p(z),&s=0\ \text{or}\ u=0,
 \end{cases}\\
 h_s(z,u)&=\ip u{q_s(z,u)-p(z)}.
 \end{aligned}
 \label{eq:hsdefs}
\end{equation}
In \eqref{eq:smoothdefs} and \eqref{eq:hsdefs}, argmin and argmax denote the unique minimizer and maximizer, respectively; their existence and uniqueness are proved in Lemma~\ref{lem:regularity}. We call $q_s(z,u)$ the support point of $C_s(z)$ in the direction $u$. By construction, $h_s(z,\cdot)$ is the support function of $C_s(z)-p(z)$ in the sense of Lemma~\ref{lem:geometry}. Gradients are column vectors, and $D_zq_s=(\partial(q_s)_i/\partial z_j)_{i,j=1}^k$ is the Jacobian matrix of $q_s$ with respect to $z$. Let $Z,U$ be independent standard Gaussian vectors in $\R^k$. A subscript on $\E$ indicates the variables over which the expectation is taken, the remaining arguments being held fixed.

\begin{lemma}\label{lem:support}
For every $\theta>0$, all terms below are finite, and we have
\begin{equation}
 \E h_\theta(Z,U)^2
 =2\int_0^\theta\E\ip Z{q_s(Z,U)} ds,
 \label{eq:supportidentity}
\end{equation}
and we also have
\begin{subequations}\label{eq:smoothbounds}
\begin{align}
 \E\diam(C_\theta(Z))^2
 &\le8\int_0^\theta\E\ip Z{q_s(Z,U)} ds,\\
 \E\max_{t\in C_\theta(Z)}\norm{t-p(Z)}^2
 &\le4\int_0^\theta\E\ip Z{q_s(Z,U)} ds.
\end{align}
\end{subequations}
\end{lemma}
\begin{proof}
By Lemma~\ref{lem:regularity} in Appendix~\ref{app:regularity}, the map $p$ is continuously differentiable with $\nabla H^*=p$. For $s>0$ and $u\ne0$, the support point $q=q_s(z,u)$, the vector $a=\nabla H(q)-z$, and a multiplier $\lambda=\lambda_s(z,u)$ satisfy
\begin{equation}
 u=\lambda a,\qquad H(q)-\ip zq+H^*(z)=s,\qquad\lambda>0,
 \label{eq:kkt}
\end{equation}
and $q_s$ and $\lambda$ are continuously differentiable on $\{s>0,u\ne0\}$.

Differentiating the second equation of \eqref{eq:kkt} with respect to $s$ and $z$, and writing $\partial_s$ for the partial derivative in $s$, we obtain
\[
 a^\top\partial_s q_s=1,\qquad
 a^\top D_zq_s=(q-p(z))^\top.
\]
Since $p(z)\in C_s(z)$, we have $h_s\ge0$. Using $u=\lambda a$ from \eqref{eq:kkt} and the last display, we obtain
\begin{align}
 \partial_s h_s
 &=u^\top\partial_sq_s=\lambda a^\top\partial_sq_s=\lambda,\notag\\
 u^\top D_zq_su
 &=\lambda a^\top D_zq_su
 =\lambda(q-p(z))^\top u=\lambda h_s,\notag\\
 \partial_s(h_s^2)
 &=2h_s\partial_sh_s=2u^\top D_zq_su\ge0.\label{eq:keyderivative}
\end{align}

For $u\ne0$, we integrate \eqref{eq:keyderivative} over $[b,\theta]$ and let $b\downarrow0$, which is allowed since the integrand is nonnegative and $h_b(z,u)\to0$ by \eqref{eq:integrability}. Since both sides below vanish for $u=0$, for all $z,u\in\R^k$ we obtain
\begin{equation}
 h_\theta(z,u)^2=2\int_0^\theta u^\top D_zq_s(z,u)u ds.
 \label{eq:integralderivative}
\end{equation}

Two Gaussian integrations by parts, justified in Lemma~\ref{lem:ibp}, remove the derivative in \eqref{eq:integralderivative}. For $s>0$ and $z,u\in\R^k$, with $\psi_s(z,u)=\ip u{q_s(z,u)}$ the support function of $C_s(z)$, we have
\begin{align}
 \E_Z[u^\top D_zq_s(Z,u)u]
 &=\E_Z[\ip Zu\ip u{q_s(Z,u)}],\label{eq:firstibp}\\
 \E_U[U \psi_s(z,U)]
 &=\E_Uq_s(z,U).\label{eq:secondibp}
\end{align}
The first is an integration by parts with respect to the component of $Z$ along $u$. The second is Stein's lemma for $\psi_s(z,\cdot)$, whose gradient is $q_s(z,\cdot)$.

The integrand in \eqref{eq:integralderivative} is nonnegative, so Tonelli's theorem applies. For $0<s\le\theta$, \eqref{eq:integrability} bounds the absolute value of the right-hand integrand in \eqref{eq:firstibp} by
\[
 \norm Z\norm U^2
 \left[\mu^{-1}\norm{Z-\nabla H(0)}+\sqrt{2\theta/\mu}\right],
\]
which is integrable in $(Z,U,s)$. Fubini's theorem, \eqref{eq:firstibp}, and \eqref{eq:secondibp} now yield
\[
 \E h_\theta(Z,U)^2
 =2\int_0^\theta\E[\ip ZU\ip U{q_s(Z,U)}] ds
 =2\int_0^\theta\E\ip Z{q_s(Z,U)} ds.
\]
Applying Lemma~\ref{lem:geometry} to $B=C_\theta(z)-p(z)$, whose support function is $h_\theta(z,\cdot)$, then setting $z=Z$ and taking expectations, we obtain \eqref{eq:smoothbounds}; the diameter and the maximal distance are Borel functions of $Z$ by Lemma~\ref{lem:regularity}.
\end{proof}

\subsection{Reduction to finite index sets}
\label{sec:regularization}

\begin{proposition}\label{thm:finite}
Let $k,m\ge1$ be integers, $v_1,\ldots,v_m\in\R^k$ not necessarily distinct, $a_1,\ldots,a_m\in\R$, and $Z\sim N(0,I_k)$. For $z\in\R^k$ and $\eps\ge0$, set
\[
 F(z)=\max_i(a_i+\ip z{v_i}),\qquad
 A_\eps(z)=\{v_i:a_i+\ip z{v_i}\ge F(z)-\eps\}.
\]
Let $p(z)$ be $v_i$ for the smallest index $i$ maximizing $a_i+\ip z{v_i}$. Then $p$ is Borel and satisfies $A_0(Z)=\{p(Z)\}$ almost surely, and, for every $\eps\ge0$, we have
\begin{align}
 \E\diam(A_\eps(Z))^2
 &\le8\eps \E\max_i\ip Z{v_i},\\
 \E\max_{x\in A_\eps(Z)}\norm{x-p(Z)}^2
 &\le4\eps \E\max_i\ip Z{v_i}.
\end{align}
If, in addition, $U\sim N(0,I_k)$ is independent of $Z$, then, for every $\eps\ge0$, we have
\begin{equation}
 \E\Bigl(\max_{x\in A_\eps(Z)}\ip U{x-p(Z)}\Bigr)^2
 \le2\eps \E\max_i\ip Z{v_i}.
 \label{eq:finitesharp}
\end{equation}
\end{proposition}

The proof is in Appendix~\ref{app:finite}. It applies Lemma~\ref{lem:support} to the conjugate $H=F_N^*$ of the regularized function $F_N(z)=N^{-1}\log\sum_{i=1}^m\exp\{N(a_i+\ip z{v_i})\}+\norm z^2/(2N)$ and lets $N\to\infty$. For fixed $z$ and large $N$, the sets $C_s(z)$ defined by this $H$ contain $A_\eps(z)$ when $s>\eps$, and they lie within a distance of order $N^{-1/2}$ of the convex hull of the $v_i$. Hence the inner product $\ip Z{q_s(Z,U)}$ in \eqref{eq:supportidentity} is at most $\max_i\ip Z{v_i}$ plus an error that vanishes in mean as $N\to\infty$.

Theorem~\ref{thm:process} is deduced from Proposition~\ref{thm:finite} in Appendix~\ref{sec:process}. For a finite index set, the vector $(G(t))_{t\in T}$ has the law of $(\ip Z{v_t})_{t\in T}$ for suitable vectors $v_t$, and the proposition applies with $\eps/\sigma$ in place of $\eps$. A compact index set is approximated by finite nets, and Fatou's lemma carries the bounds to the limit.

\section{Optimality and dependence on the drift}
\label{sec:sharpness}

This section examines how sharp the bounds of Theorem~\ref{thm:process} are. Section~\ref{sec:optimality} states that the constants $8$, $4$, and $2$ cannot be improved, with the examples in which they are approached; Remark~\ref{rem:balls} compares the bounds with the limits as $\eps\downarrow0$ for Euclidean balls. Section~\ref{sec:drift} concerns a given drift: Proposition~\ref{prop:adaptive} states bounds that adapt to the drift and that replace $\E\sup_TG$ by $\E\max_{A_0}G$ in the limit $\eps\downarrow0$, and Example~\ref{ex:chernoff} applies them to Brownian motion with parabolic drift. Remark~\ref{rem:assumptions} discusses the assumptions of Theorem~\ref{thm:process}. The proofs are in Appendices~\ref{app:optimality} and~\ref{app:adaptive}.

\subsection{Optimality of the constants}
\label{sec:optimality}

\begin{proposition}\label{prop:sharpness}
Let $\eps>0$. If the constant $8$ in \eqref{eq:processbound}, $4$ in \eqref{eq:processradius}, or $2$ in \eqref{eq:processsharp} is replaced by a smaller one, then the inequality fails for some compact $T$ and some $G$, with $f=0$ and $\sigma=1$.
\end{proposition}

The proof is in Appendix~\ref{app:optimality}. With $f=0$ and $\sigma=1$, the ratio of the two sides of \eqref{eq:processbound} tends to $1$ as $\eps\downarrow0$ when $T$ is the unit disk of $\R^2$ and $G(t)=\ip Zt$ with $Z\sim N(0,I_2)$. The same holds for \eqref{eq:processradius} and \eqref{eq:processsharp} when $T=[0,1]$ and $G(t)=Zt$ with scalar $Z\sim N(0,1)$. Replacing $T$ by $cT$ and letting $c\to\infty$ gives the assertion for a fixed $\eps$.

\begin{remark}\label{rem:balls}
For a given index set, the constants are approached in these examples only as $\eps\downarrow0$, and \eqref{eq:processbound} can be far from sharp. Let $T$ be the Euclidean unit ball of $\R^k$ with $k\ge2$, $G(t)=\ip Zt$ with $Z\sim N(0,I_k)$, $f=0$, and $\sigma=1$. Appendix~\ref{app:optimality} shows that
\begin{equation}
 \begin{aligned}
 \lim_{\eps\downarrow0}\frac{\E\diam_d(A_\eps)^2}{\eps}
 &=\frac{8}{k-1} \E\sup_TG,\\
 \lim_{\eps\downarrow0}\frac1\eps 
 \E\max_{\substack{t\in A_\eps\\p\in A_0}}d(t,p)^2
 &=\frac{2}{k-1} \E\sup_TG,\\
 \lim_{\eps\downarrow0}\frac{\E\Gamma_\eps^2}{\eps}
 &=2 \E\sup_TG.
 \end{aligned}
 \label{eq:balllimits}
\end{equation}
Thus \eqref{eq:processbound} and \eqref{eq:processradius} overestimate the first two limits by the factors $k-1$ and $2(k-1)$, whereas \eqref{eq:processsharp} is attained asymptotically for this ball. Replacing $T$ by $b+cT$ with $b\in\R^k$ and $c>0$, or $\sigma=1$ by another noise level, only rescales $\eps$, so \eqref{eq:processsharp} is attained asymptotically for every Euclidean ball of dimension $k\ge2$ and every $\sigma>0$; for $k=1$ this follows from the segment example in the proof of Proposition~\ref{prop:sharpness}. For balls, the loss in \eqref{eq:processbound} therefore arises entirely in Lemma~\ref{lem:geometry}.
\end{remark}

\subsection{Dependence on the drift}
\label{sec:drift}

For a particular drift, the bounds of Theorem~\ref{thm:process} may be far from sharp. Take $T=[0,1]$, $G(t)=Zt$ with scalar $Z\sim N(0,1)$, $\sigma=1$, and $f(t)=-\kappa t$ with $\kappa\ge0$. Then we have $X(t)=(Z-\kappa)t$, so the maximal distance from the maximizer equals $\min\{1,\eps/\abs{Z-\kappa}\}$ outside $\{Z=\kappa\}$. The dominated convergence argument of Step~2 in the proof of Proposition~\ref{prop:sharpness}, with $\ph(\eps u)$ replaced by $\ph(\kappa\pm\eps u)$ on the two half-lines, gives
\[
 \lim_{\eps\downarrow0}\frac1\eps 
 \E\max_{\substack{t\in A_\eps\\p\in A_0}}\abs{t-p}^2
 =\lim_{\eps\downarrow0}\frac1\eps \E\min\{1,\eps/\abs{Z-\kappa}\}^2
 =4\ph(\kappa),
\]
whereas \eqref{eq:processradius} bounds this limit by $4\ph(0)$ for every $\kappa$.

The proof of Theorem~\ref{thm:process} yields a family of bounds that adapts to the drift. The random variable $\max_{t\in A_0}G(t)$ below is measurable by the upper semicontinuity argument of Appendix~\ref{sec:process}, applied to the map $g\mapsto\max\{g(t):t\in A_0(f+\sigma g)\}$ on $C(T)$, with $A_0(\cdot)$ as in the proof of Theorem~\ref{thm:process}; when $d$ is a metric it equals $G(p)$ for the almost surely unique maximizer $p$.

\begin{proposition}\label{prop:adaptive}
In the setting of Theorem~\ref{thm:process}, put
\[
 \Delta(\lambda)=\E\max_T\bigl(f+(1+\lambda)\sigma G\bigr)-\E\max_T\bigl(f+\sigma G\bigr),
 \qquad\lambda>0.
\]
Then $\Delta(\lambda)/\lambda$ is nondecreasing in $\lambda$, with $\Delta(\lambda)/\lambda\to\sigma\E\max_{t\in A_0}G(t)$ as $\lambda\downarrow0$ and $\Delta(\lambda)/\lambda\to\sigma\E\sup_TG$ as $\lambda\to\infty$. For every $\eps\ge0$ and $\lambda>0$, with $\Gamma_\eps$ as in Theorem~\ref{thm:process}, we have
\begin{align}
 \E\diam_d(A_\eps)^2
 &\le\frac{8\eps \Delta(\lambda)+4\eps^2}{\lambda\sigma^2},\label{eq:adaptivediam}\\
 \E\max_{\substack{t\in A_\eps\\p\in A_0}}d(t,p)^2
 &\le\frac{4\eps \Delta(\lambda)+2\eps^2}{\lambda\sigma^2},\label{eq:adaptiveradius}\\
 \E\Gamma_\eps^2
 &\le\frac{2\eps \Delta(\lambda)+\eps^2}{\lambda\sigma^2}.\label{eq:adaptivesharp}
\end{align}
Consequently, we have
\begin{equation}
 \begin{aligned}
 \limsup_{\eps\downarrow0}\frac{\E\diam_d(A_\eps)^2}{\eps}
 &\le\frac8\sigma\E\max_{t\in A_0}G(t),\\
 \limsup_{\eps\downarrow0}\frac1\eps 
 \E\max_{\substack{t\in A_\eps\\p\in A_0}}d(t,p)^2
 &\le\frac4\sigma\E\max_{t\in A_0}G(t),\\
 \limsup_{\eps\downarrow0}\frac{\E\Gamma_\eps^2}{\eps}
 &\le\frac2\sigma\E\max_{t\in A_0}G(t).
 \end{aligned}
 \label{eq:adaptivelimit}
\end{equation}
\end{proposition}

Proposition~\ref{prop:adaptive} is proved in Appendix~\ref{app:adaptive}. Letting $\lambda\to\infty$ in \eqref{eq:adaptivediam}--\eqref{eq:adaptivesharp} recovers \eqref{eq:processsharp} and the first terms of the minima in \eqref{eq:processbound} and \eqref{eq:processradius}. Moreover, we have $0\le\E\max_{A_0}G\le\E\sup_TG$, since $\max_{A_0}G\ge(\max_TX-\max_Tf)/\sigma$, whose expectation is nonnegative, and $A_0\subset T$; the asymptotic constants in \eqref{eq:adaptivelimit} are therefore at most those of Theorem~\ref{thm:process}, with equality when $f$ is constant. In the example above, we have $\E\max_{A_0}G=\E[Z\1_{\{Z>\kappa\}}]=\ph(\kappa)$, so the second inequality in \eqref{eq:adaptivelimit} is an equality for every $\kappa\ge0$, and so is the third, since $\E[\Gamma_\eps^2\mid Z]$ is half the squared maximal distance for $\widetilde G(t)=Ut$; the first is not, since there $\diam(A_\eps)$ equals the maximal distance, which gives $\eps^{-1}\E\diam(A_\eps)^2\to4\ph(\kappa)$.

\begin{example}\label{ex:chernoff}
Let $W$ be a two-sided standard Brownian motion, and let $L>0$, $T=[-L,L]$, $G(t)=W(t)$, $f(t)=-t^2$, and $\sigma=1$. Then we have $d(s,t)=\abs{s-t}^{1/2}$, so the maximizer $p$ of $X$ is almost surely unique. The reflection principle, applied to the two independent halves of $W$, gives $\E\sup_TG=2\sqrt{L/\pi}$, so \eqref{eq:processradius} yields $\E\max_{t\in A_\eps}\abs{t-p}\le8\eps\sqrt{L/\pi}$, which grows with $L$. The bound \eqref{eq:adaptiveradius} does not. Since $t^2\ge2\abs t-1$ and $\sup_{t\ge0}(W(t)-t)$ is exponentially distributed with mean $1/2$, we have $\E\sup_{t\in\R}(2W(t)-t^2)\le3$, which gives $\Delta(1)\le3$ because $\max_TX\ge X(0)=0$. Therefore, for every $L>0$ and $\eps\ge0$, \eqref{eq:adaptiveradius} with $\lambda=1$ gives
\[
 \E\max_{t\in A_\eps}\abs{t-p}\le12\eps+2\eps^2.
\]
This bound passes to the whole line. Since $W(t)-t^2\to-\infty$ as $\abs t\to\infty$ almost surely, for almost every path the set $A_\eps$ and the point $p$ no longer change once $L$ is large enough: they are then the set of $\eps$-near-maximizers of $W(t)-t^2$ over $\R$ and its maximizer, whose law is Chernoff's distribution~\cite{Chernoff1964,Groeneboom1989}. By Fatou's lemma, the display holds for these as well. The uniqueness of the maximizer is for the fixed drift: the largest maximizer of $W(t)-(t-a)^2$ over $\R$ is a jump process in $a\in\R$~\cite{Groeneboom1989}, and at each of its jump times there are two maximizers.
\end{example}

\begin{remark}\label{rem:assumptions}
Theorem~\ref{thm:process} extends to a random $\widehat\eps\ge0$ that is almost surely bounded by a deterministic constant, which then replaces $\eps$. By conditioning, it also extends to a random $\widehat\eps$ and a random drift that are independent of $(G,\widetilde G)$, with $\E\widehat\eps$ in place of $\eps$; by the argument of Step~1 in the proof of Theorem~\ref{thm:process}, the random variables in the theorem are jointly Borel in $\eps$, $f$, and the paths. Neither $\widehat\eps$ nor the drift may depend on $G$. For $G(t)=Zt$ on $[0,1]$ with scalar $Z\sim N(0,1)$ and $\sigma=1$, the drift $f(t)=-Zt$ gives $X\equiv0$ and $A_0=[0,1]$; and for $f=0$, the choice $\widehat\eps=\abs Z\1_{\{\abs Z\le\eta\}}$ gives $A_{\widehat\eps}=[0,1]$ on an event of probability $2\ph(0)\eta+o(\eta)$ as $\eta\downarrow0$, whereas we have $\E\widehat\eps=\ph(0)\eta^2+o(\eta^2)$. The bounds also do not pass to metrics other than $d$ without a comparison: for $G(t)=Zt^2$ on $[-1,1]$ and $f=0$, the points $\pm1$ are at canonical distance zero and both maximize $X$ on $\{Z>0\}$.
\end{remark}

\section*{Acknowledgments}
This work was supported by JSPS KAKENHI Grant Number JP25K23133.

The author used ChatGPT 6.1-Sol (OpenAI) and Claude Opus 5.5 (Anthropic) to assist in checking the proofs and in editing the English. All mathematical statements and proofs were verified by the author, who takes full responsibility for the content of this article.

\newpage

\appendix

\section{Regularity of the support points}
\label{app:regularity}

This appendix proves the regularity and the two integrations by parts used in the proof of Lemma~\ref{lem:support}, in the setting and notation of Section~\ref{sec:support}.

\begin{lemma}\label{lem:regularity}
(i) For every $z\in\R^k$, the function $t\mapsto H(t)-\ip zt$ has a unique minimizer $p(z)$, which satisfies $\nabla H(p(z))=z$. The map $p$ is continuously differentiable, and we have $\nabla H^*=p$ and $C_0(z)=\{p(z)\}$. (ii) For $s>0$ and $z\in\R^k$, the set $C_s(z)$ is compact and strictly convex. For $u\ne0$, the linear functional $\ip u\cdot$ attains its maximum over $C_s(z)$ at a unique point $q=q_s(z,u)$, and there is $\lambda>0$ such that $q$, $a=\nabla H(q)-z$, and $\lambda$ satisfy \eqref{eq:kkt}. As functions of $(s,z,u)$, $q_s$ and $\lambda$ are continuously differentiable on $\{s>0,u\ne0\}$. (iii) For $s\ge0$ and $z,u\in\R^k$, we have
\begin{equation}
 \max_{t\in C_s(z)}\norm t\le\mu^{-1}\norm{z-\nabla H(0)}+\sqrt{2s/\mu},
 \qquad 0\le h_s(z,u)^2\le(2s/\mu)\norm u^2.
 \label{eq:integrability}
\end{equation}
(iv) The map $(s,z,u)\mapsto q_s(z,u)$ is Borel. For every $\theta>0$, the diameter of $C_\theta(z)$ and the maximal distance from $p(z)$ to a point of $C_\theta(z)$ are upper semicontinuous, hence Borel, functions of $z$.
\end{lemma}
\begin{proof}
Strong convexity makes $H(t)-\ip zt$ coercive, with a unique minimizer $p(z)$ satisfying $\nabla H(p(z))=z$. By the inverse function theorem $p$ is continuously differentiable, and differentiating $H^*(z)=\ip z{p(z)}-H(p(z))$ gives $\nabla H^*=p$. Moreover, by the strong convexity of $H$, $p$ and the Fenchel--Young gap satisfy
\begin{align}
 \mu\norm{p(z)}^2
 &\le\ip{z-\nabla H(0)}{p(z)}
 \le\norm{z-\nabla H(0)}\norm{p(z)},\label{eq:pbound}\\
 H(t)-\ip zt+H^*(z)
 &=H(t)-H(p(z))-\ip{\nabla H(p(z))}{t-p(z)}\notag\\
 &\ge\frac\mu2\norm{t-p(z)}^2.\label{eq:stronggap}
\end{align}
In particular, we have $C_0(z)=\{p(z)\}$.

For $s>0$, $C_s(z)$ is compact, has $p(z)$ in its interior, and is strictly convex: the midpoint of distinct $x,y\in C_s(z)$ has gap at most $s-\mu\norm{x-y}^2/8<s$. Hence, for $u\ne0$, the linear functional $\ip u\cdot$ attains its maximum over $C_s(z)$ at a unique point $q=q_s(z,u)$ on the boundary; since $p(z)$ is interior, we have $q\ne p(z)$, and since $\nabla H(t)=z$ only at $t=p(z)$, this gives $a=\nabla H(q)-z\ne0$. Hence the boundary is a continuously differentiable hypersurface near $q$ with tangent space $a^\perp$, and the Lagrange multiplier rule gives $u=\lambda a$ for some $\lambda\in\R$; since $\ip u{q-p(z)}>0$ and $\ip a{q-p(z)}\ge\mu\norm{q-p(z)}^2>0$, we have $\lambda>0$. This gives \eqref{eq:kkt}.

The Jacobian in $(q,\lambda)$ of the residuals $\lambda(\nabla H(q)-z)-u$ and $H(q)-\ip zq+H^*(z)-s$ is $\left(\begin{smallmatrix}\lambda\nabla^2H(q)&a\\a^\top&0\end{smallmatrix}\right)$. If it maps $(\xi,\xi_0)\in\R^k\times\R$ to zero, then we have $\lambda\nabla^2H(q)\xi+\xi_0a=0$ and $a^\top\xi=0$, which give $\lambda\xi^\top\nabla^2H(q)\xi=0$, $\xi=0$, and $\xi_0=0$; hence the matrix is invertible. The implicit function theorem therefore yields, near each $(s,z,u)$ with $s>0,u\ne0$, a continuously differentiable branch of solutions of \eqref{eq:kkt} along which $\lambda$ stays positive. For any solution $(q,\lambda)$ of \eqref{eq:kkt} with $\lambda>0$, the point $q$ is the support point $q_s(z,u)$: for $t\in C_s(z)$, convexity of the gap gives $s\ge H(t)-\ip zt+H^*(z)\ge s+\ip a{t-q}=s+\lambda^{-1}\ip u{t-q}$, which yields $\ip ut\le\ip uq$. The branch therefore coincides with $(q_s(z,u),\lambda)$, and $q_s$ and $\lambda$ are continuously differentiable on $\{s>0,u\ne0\}$.

The bounds \eqref{eq:integrability} follow from \eqref{eq:pbound} and \eqref{eq:stronggap}. The map $(s,z,u)\mapsto q_s(z,u)$ is Borel, being continuous on the open set $\{s>0,u\ne0\}$ and equal to the continuous map $p(z)$ on the closed complement. The diameter of $C_\theta(z)$ and the maximal distance from $p(z)$ are upper semicontinuous, hence Borel, functions of $z$: if $z_n\to z$, \eqref{eq:integrability} places the sets $C_\theta(z_n)$ in a common compact ball, and by continuity of $H$, $H^*$, and $p$ the limit of any convergent sequence of maximizing pairs or points belongs to $C_\theta(z)$.
\end{proof}

\begin{lemma}\label{lem:ibp}
For every $s>0$, the identity \eqref{eq:firstibp} holds for every $u\in\R^k$ and the identity \eqref{eq:secondibp} holds for every $z\in\R^k$, and the expectations in them are finite.
\end{lemma}
\begin{proof}
For $u=0$ both sides of \eqref{eq:firstibp} vanish. Fix $s>0$ and $u\ne0$, put $v=u/\norm u$, and write $Z=Z_\perp+Sv$, where $S\sim N(0,1)$ is independent of the orthogonal projection $Z_\perp$ of $Z$ onto $v^\perp$. For $y\in v^\perp$, set $g_y(r)=\ip u{q_s(y+rv,u)}$; then we have $g_y'(r)=\norm u^{-1}u^\top D_zq_s(y+rv,u)u\ge0$ by \eqref{eq:keyderivative} and $\abs{g_y(r)}\le\norm u[\mu^{-1}(\norm y+\abs r+\norm{\nabla H(0)})+\sqrt{2s/\mu}]$ by \eqref{eq:integrability}. We integrate by parts on $[-R,R]$ and let $R\to\infty$: the boundary term $[g_y(r)\ph(r)]_{-R}^R$ vanishes by the linear growth bound, monotone convergence applies to $\int g_y'\ph$, and dominated convergence to $\int rg_y\ph$; this gives $\E_Sg_y'(S)=\E_S[Sg_y(S)]$. The same growth bound justifies averaging over $Z_\perp$, and multiplying by $\norm u$ gives \eqref{eq:firstibp}.

For fixed $z$ and $s>0$, the support function $\psi_s(z,u)=\ip u{q_s(z,u)}=h_s(z,u)+\ip u{p(z)}$ of $C_s(z)$ is Lipschitz on $\R^k$, and for $u\ne0$ it satisfies $\ip{u'-u}{q_s(z,u)}\le\psi_s(z,u')-\psi_s(z,u)\le\ip{u'-u}{q_s(z,u')}$ for all $u'\in\R^k$, so continuity of $q_s(z,\cdot)$ at $u\ne0$ gives $\nabla_u\psi_s(z,u)=q_s(z,u)$. On every line parallel to a coordinate axis, $\psi_s(z,\cdot)$ is therefore absolutely continuous, with bounded derivative equal to the corresponding coordinate of $q_s(z,\cdot)$ except possibly at the origin. Stein's lemma~\cite{Stein1981}, applied to each coordinate of $U$ conditionally on the other coordinates, gives \eqref{eq:secondibp}.
\end{proof}

\section{Proof for finite index sets}
\label{app:finite}

\begin{proof}[Proof of Proposition~\ref{thm:finite}]
By definition, we have
\[
 p(z)=\sum_{i=1}^m v_i 
 \1_{\{a_i+\ip z{v_i}=F(z),\
             a_j+\ip z{v_j}<F(z)\ \text{for all }j<i\}},
\]
so $p$ is Borel. If $v_i\ne v_j$, then $a_i-a_j+\ip Z{v_i-v_j}$ is Gaussian with variance $\norm{v_i-v_j}^2>0$. This gives
\[
 \Pbb\{\abs{A_0(Z)}>1\}
 \le\sum_{\substack{i<j\\v_i\ne v_j}}
   \Pbb\{a_i-a_j+\ip Z{v_i-v_j}=0\}=0.
\]

Subtracting $\min_j a_j$ from each $a_i$ leaves $A_\eps$ and $p$ unchanged, so we may assume $\min_i a_i=0$. Let $K=\conv\{v_1,\ldots,v_m\}$ be their convex hull. For $N=1,2,\ldots$ define
\[
 F_N(z)=\frac1N\log\sum_{i=1}^m
 \exp\{N(a_i+\ip z{v_i})\}+\frac{\norm z^2}{2N}.
\]
For $\pi_{N,i}(z)=\exp\{N(a_i+\ip z{v_i})\}/\sum_j\exp\{N(a_j+\ip z{v_j})\}$, direct differentiation gives
\begin{equation}
 \nabla F_N(z)=\sum_i\pi_{N,i}(z)v_i+z/N.
 \label{eq:softmax-gradient}
\end{equation}
Since $\pi_{N,i}\ge0$ and $\sum_i\pi_{N,i}=1$, the bracket below is a variance under the weights $\pi_{N,i}$, so for all $z,u\in\R^k$ we have
\begin{gather*}
 u^\top\nabla^2F_N(z)u
 =N\Bigl[\sum_i\pi_{N,i}\ip u{v_i}^{2}
 -\Bigl(\sum_i\pi_{N,i}\ip u{v_i}\Bigr)^2\Bigr]+\frac{\norm u^2}{N},\\
 \frac{\norm u^2}{N}\le u^\top\nabla^2F_N(z)u
 \le\Bigl(N\max_i\norm{v_i}^2+\frac1N\Bigr)\norm u^2.
\end{gather*}
For fixed $N$ and $x\in\R^k$, the function $F_N(z)-\ip xz$ is strictly convex and, being bounded below by $a_1+\ip z{v_1-x}+\norm z^2/(2N)$, coercive; its unique minimizer $y(x)$ satisfies $\nabla F_N(y(x))=x$. The inverse function theorem and differentiation of $F_N^*(x)=\ip x{y(x)}-F_N(y(x))$ give
\begin{gather*}
 \nabla F_N^*(x)=y(x),\qquad
 \nabla^2F_N^*(x)=\nabla^2F_N(y(x))^{-1},\\
 u^\top\nabla^2F_N^*(x)u
 \ge\frac{\norm u^2}{N\max_i\norm{v_i}^2+N^{-1}},
 \qquad x,u\in\R^k.
\end{gather*}
Thus $F_N^*$ is smooth and strongly convex. By definition, we have $F_N^*(x)\ge\ip zx-F_N(z)$ for all $x,z\in\R^k$, with equality at $x=\nabla F_N(z)$; this gives $F_N^{**}=F_N$. Lemma~\ref{lem:support} therefore applies with $H=F_N^*$. The point denoted $p(z)$ in that lemma is then $\nabla F_N(z)$, whereas $p(z)$ keeps its meaning from the statement of the proposition.

The elementary log-sum-exp bounds give
\begin{gather}
 F(z)+\frac{\norm z^2}{2N}
 \le F_N(z)\le F(z)+\frac{\log m}{N}+\frac{\norm z^2}{2N},
 \qquad F_N^*(v_i)\le-a_i.
 \label{eq:reg-conjugate}
\end{gather}
The last inequality follows by taking the supremum of $\ip y{v_i}-F_N(y)\le-a_i-\norm y^2/(2N)$ over $y$. Writing $\dist(x,K)=\min_{\bar x\in K}\norm{x-\bar x}$, we also have
\begin{equation}
 F_N^*(x)\ge-\max_i a_i-\frac{\log m}{N}
              +\frac N2\dist(x,K)^2.
 \label{eq:reg-lower}
\end{equation}
Indeed, let $\bar x\in K$ be nearest to $x$ and put $y=N(x-\bar x)$. The variational inequality for the projection onto the closed convex set $K$ gives $\max_i\ip y{v_i}=\max_{x'\in K}\ip y{x'}=\ip y{\bar x}$, so \eqref{eq:reg-conjugate} gives $F_N(y)\le\max_ia_i+\ip y{\bar x}+(\log m)/N+\norm y^2/(2N)$, and substituting this into $F_N^*(x)\ge\ip yx-F_N(y)$ proves \eqref{eq:reg-lower}.

For $s\ge0$ and $z,u\in\R^k$, let $C_{N,s}(z)$ and $h_{N,s}(z,u)$ denote the set $C_s(z)$ of \eqref{eq:smoothdefs} and the function $h_s(z,u)$ of \eqref{eq:hsdefs} for $H=F_N^*$, so that we have
\begin{gather*}
 C_{N,s}(z)=\{x\in\R^k:F_N^*(x)-\ip zx+F_N(z)\le s\},\\
 h_{N,s}(z,u)=\max_{x\in C_{N,s}(z)}\ip u{x-\nabla F_N(z)}.
\end{gather*}
For $x\in C_{N,s}(z)$, let $\bar x\in K$ be nearest to $x$. Since $F_N(z)\ge F(z)\ge\max_i\ip z{v_i}\ge\ip z{\bar x}$, the second inequality because $a_i\ge0$, \eqref{eq:reg-lower} yields
\[
 s\ge F_N^*(x)-\ip zx+F_N(z)
 \ge\frac N2\dist(x,K)^2-\norm z\dist(x,K)-\max_i a_i-\frac{\log m}{N}.
\]
Solving the quadratic inequality and using the subadditivity of the square root gives
\begin{align*}
 \dist(x,K)&\le
 \frac{\norm z+\sqrt{\norm z^2+2N\{s+\max_i a_i+(\log m)/N\}}}{N}\\
 &\le\frac{2\norm z}{N}
       +\sqrt{\frac{2\{s+\max_i a_i+(\log m)/N\}}{N}}.
\end{align*}
Thus every $x\in C_{N,s}(z)$ satisfies
\[
 \ip zx\le\max_i\ip z{v_i}+\frac{2\norm z^2}{N}
       +\norm z\sqrt{\frac{2\{s+\max_i a_i+(\log m)/N\}}{N}}.
\]
For fixed $s>0$, the maximum $\max_{x\in C_{N,s}(z)}\ip zx$ equals $\ip z{q_s(z,z)}$, with the support point $q_s$ of \eqref{eq:hsdefs} for $H=F_N^*$. It is therefore a Borel function of $z$; its absolute value is at most $\norm z\max_{x\in C_{N,s}(z)}\norm x$, which is integrable at $z=Z$ by \eqref{eq:integrability}. For fixed $\theta>0$, taking the maximum over $x\in C_{N,s}(z)$, setting $z=Z$, and taking expectations, we obtain
\begin{equation}
 \begin{aligned}
 \sup_{0<s\le\theta}\E\max_{x\in C_{N,s}(Z)}\ip Zx
 &\le\E\max_i\ip Z{v_i}+\frac{2k}{N}\\
 &\quad+\E\norm Z\sqrt{\frac{2\{\theta+\max_i a_i+(\log m)/N\}}{N}}.
 \end{aligned}
 \label{eq:widthbudget}
\end{equation}
The right-hand side tends to $\E\max_i\ip Z{v_i}$ as $N\to\infty$. For $H=F_N^*$, the support points in \eqref{eq:supportidentity} satisfy $q_s(Z,U)\in C_{N,s}(Z)$. Hence Lemma~\ref{lem:support} and \eqref{eq:widthbudget} give
\begin{subequations}\label{eq:regbounds}
\begin{align}
 \limsup_N\E h_{N,\theta}(Z,U)^2
 &\le2\theta \E\max_i\ip Z{v_i},\\
 \limsup_N\E\diam(C_{N,\theta}(Z))^2
 &\le8\theta \E\max_i\ip Z{v_i},\\
 \limsup_N\E\max_{x\in C_{N,\theta}(Z)}
       \norm{x-\nabla F_N(Z)}^2
 &\le4\theta \E\max_i\ip Z{v_i}.
\end{align}
\end{subequations}

For every index $i$ satisfying $a_i+\ip z{v_i}\ge F(z)-\eps$, \eqref{eq:reg-conjugate} gives
\[
 F_N^*(v_i)-\ip z{v_i}+F_N(z)
 \le-a_i-\ip z{v_i}+F(z)+\frac{\log m}{N}+\frac{\norm z^2}{2N}
 \le\eps+\frac{\log m}{N}+\frac{\norm z^2}{2N}.
\]
Every point of $A_\eps(z)$ is of the form $v_i$ for such an index $i$. Hence, for fixed $z$ and $\delta>0$, we have $A_\eps(z)\subset C_{N,\eps+\delta}(z)$ for all sufficiently large $N$. If $A_0(z)=\{p(z)\}$, then we have $F(z)-a_j-\ip z{v_j}>0$ whenever $v_j\ne p(z)$. As $N\to\infty$, this gives
\begin{gather*}
 \pi_{N,j}(z)
 \le\exp\{-N(F(z)-a_j-\ip z{v_j})\}\longrightarrow0
 \qquad(v_j\ne p(z)),\\
 \nabla F_N(z)-p(z)
 =\sum_{j:v_j\ne p(z)}\pi_{N,j}(z)(v_j-p(z))+z/N
 \longrightarrow0.
\end{gather*}

The inclusion for large $N$ gives $\diam(A_\eps(z))^2\le\liminf_N\diam(C_{N,\eps+\delta}(z))^2$. If $A_0(z)=\{p(z)\}$, choose $x\in A_\eps(z)$ maximizing $\norm{x-p(z)}$; then we have $\max_{x'\in A_\eps(z)}\norm{x'-p(z)}^2=\lim_N\norm{x-\nabla F_N(z)}^2\le\liminf_N\max_{x'\in C_{N,\eps+\delta}(z)}\norm{x'-\nabla F_N(z)}^2$. Moreover, in the same case, for every $u\in\R^k$ and all large $N$, the inclusion gives $\max_{x\in A_\eps(z)}\ip u{x-p(z)}\le h_{N,\eps+\delta}(z,u)+\ip u{\nabla F_N(z)-p(z)}$, where the last term tends to zero; both the left-hand side and $h_{N,\eps+\delta}(z,u)$ are nonnegative, since $p(z)\in A_\eps(z)$ and $\nabla F_N(z)\in C_{N,\eps+\delta}(z)$, so the square of the left-hand side is at most $\liminf_Nh_{N,\eps+\delta}(z,u)^2$. The maxima over $A_\eps(z)$ are Borel functions of $z$, the candidates being finite in number and selected by Borel conditions. Setting $z=Z$ and using that $A_0(Z)=\{p(Z)\}$ almost surely, we obtain from Fatou's lemma and \eqref{eq:regbounds} the three asserted bounds with $\eps+\delta$ in place of $\eps$; letting $\delta\downarrow0$ completes the proof.
\end{proof}

The regularized bounds have a classical counterpart for Gibbs averages. For fixed $z$, let $v,v'$ be independent random vectors (replicas), each taking the value $v_i$ with probability $\pi_{N,i}(z)$, $i=1,\ldots,m$, and write $\langle\cdot\rangle_z$ for the corresponding expectation, the Gibbs average at $z$. Since $\nabla^2F_N(z)$ is $N$ times the covariance matrix of $v$ plus $I_k/N$, we have $N\langle\norm{v-v'}^2\rangle_z=2\operatorname{tr}\nabla^2F_N(z)-2k/N$, and Gaussian integration by parts, $\E\operatorname{tr}\nabla^2F_N(Z)=\E\ip Z{\nabla F_N(Z)}$, together with \eqref{eq:softmax-gradient} gives
\[
 N \E\langle\norm{v-v'}^2\rangle_Z
 =2\E\ip Z{\nabla F_N(Z)}-\frac{2k}{N}
 \le2\E\max_i\ip Z{v_i}.
\]
This is \eqref{eq:finitesharp} with the Gibbs average of $\ip U{v-v'}^2$ in place of the squared maximum of $\ip U{x-p(Z)}$ over $A_\eps(Z)$, with $1/N$ in place of $\eps$, and with the same constant $2$. The maximum over the near-maximizers is what requires the identity of Lemma~\ref{lem:support}: the Gibbs average alone does not control this maximum, since the weight $\pi_{N,i}(z)$ of a distant near-maximizer $v_i$ can be small.

\section{Proof of Theorem~\ref{thm:process}}
\label{sec:process}

\begin{proof}[Proof of Theorem~\ref{thm:process}]
\emph{Step 1 (Continuity of $d$ and measurability).} If $t_n\to t$ in $T$, the increments $G(t_n)-G(t)$ are centered Gaussian and tend to zero almost surely, so they tend to zero in $L^2$, that is, $d(t_n,t)\to0$; by the triangle inequality, $d$ is continuous on $T\times T$. For $g\in C(T)$ and $\theta,\theta'\ge0$, let $A_\theta(g)=\{t\in T:g(t)\ge\max_Tg-\theta\}$. If $g_n\to g$ uniformly, every subsequential limit of a sequence $t_n\in A_\theta(g_n)$ lies in $A_\theta(g)$. By compactness and the continuity of $d$, the maps $g\mapsto\max\{d(s,t):s\in A_\theta(g), t\in A_{\theta'}(g)\}$ and $(g,\tilde g)\mapsto\max\{\tilde g(t):t\in A_\theta(g)\}$ are therefore upper semicontinuous, hence Borel, on $C(T)$ and on $C(T)\times C(T)$, respectively. Since $C(T)$ is separable, $(G,\widetilde G)$ is a Borel random element of $C(T)\times C(T)$, and all random variables in Theorem~\ref{thm:process} are measurable.

\emph{Step 2 (Finite index sets).} For finite $T=\{t_1,\ldots,t_m\}$, let $v_i$ be row $i$ of the positive semidefinite square root of the covariance matrix of $(G(t_1),\ldots,G(t_m))$. If $Z\sim N(0,I_m)$, then, with $\stackrel{d}{=}$ denoting equality in law, we have
\begin{gather*}
 (G(t_i))_{i=1}^m\ \stackrel{d}{=}\ (\ip Z{v_i})_{i=1}^m,\\
 d(t_i,t_j)^2=\E\ip Z{v_i-v_j}^{2}=\norm{v_i-v_j}^2,
 \qquad \E\max_i\ip Z{v_i}=\E\max_iG(t_i).
\end{gather*}
Apply Proposition~\ref{thm:finite} with $k=m$, $a_i=f(t_i)/\sigma$, and $\eps/\sigma$ in place of $\eps$. For $z\in\R^m$ and $\theta\ge0$, we have $\{v_i:f(t_i)+\sigma\ip z{v_i}\ge\max_j(f(t_j)+\sigma\ip z{v_j})-\theta\}=A_{\theta/\sigma}(z)$, and $d(t_i,t_j)=\norm{v_i-v_j}$. With $\theta=\eps$ and $\theta=0$, this gives \eqref{eq:processbound} and \eqref{eq:processradius} for finite $T$; repeated vectors are allowed in Proposition~\ref{thm:finite}. Representing $\widetilde G$ in the same way by a standard Gaussian vector $U$ independent of $Z$, and using $A_0(Z)=\{p(Z)\}$ almost surely, we see that $\Gamma_\eps$ has the law of $\max_{x\in A_{\eps/\sigma}(Z)}\ip U{x-p(Z)}$; thus \eqref{eq:finitesharp} gives \eqref{eq:processsharp} for finite $T$.

\emph{Step 3 (Approximation by finite nets).} For general $T$, compactness gives an increasing sequence of finite $1/n$-nets $T_n\subset T$, that is, finite sets with $\min_{s\in T_n}\rho(t,s)<1/n$ for all $t\in T$; to construct them, add the centers of a finite cover of $T$ by balls of radius $1/n$ to $T_{n-1}$. For $\theta\ge0$, put $A_{n,\theta}=\{t\in T_n:X(t)\ge\max_{T_n}X-\theta\}$. Fix $\eps\ge0$ and $\delta>0$. Let $Q$ be one of $\diam_d(A_\eps)$, $\max_{t\in A_\eps,p\in A_0}d(t,p)$, and $\Gamma_\eps$, and let $Q_n$ be the same quantity with $(T_n,\eps+\delta)$ in place of $(T,\eps)$. Steps 4--6 show that $0\le Q\le\liminf_nQ_n$ almost surely. By Fatou's lemma, Step 2, and $\E\max_{T_n}G\le\E\sup_TG$, we then have $\E Q^2\le\liminf_n\E Q_n^2\le c (\eps+\delta)\E\sup_TG/\sigma$ with $c=8$, $4$, and $2$, respectively, and $\delta\downarrow0$ completes the proof; the second terms of the minima hold since $Q\le\diam_d(T)$ in the first two cases. Two facts are used repeatedly: if $t_n\in T_n$ converge to $t\in A_\eps$, then $t_n$ belongs to $A_{n,\eps+\delta}$ eventually, by continuity of $X$ at $t$; and since $\max_{T_n}X\to\max_TX$ by uniform continuity and the net property, every subsequential limit of points of $A_{n,0}$ lies in $A_0$.

\emph{Step 4 (Diameter).} Choose $s,t\in A_\eps$ attaining the $d$-diameter and $s_n,t_n\in T_n$ converging to them. Continuity of $d$ gives $\diam_d(A_\eps)^2=\lim_nd(s_n,t_n)^2$, which is at most $\liminf_n\diam_d(A_{n,\eps+\delta})^2$. At $\eps=0$ the resulting bound is $8\delta \E\sup_TG/\sigma$ for every $\delta>0$, so we have $\diam_d(A_0)=0$ almost surely.

\emph{Step 5 (Distance from the maximizers).} By Step 4, we have $\diam_d(A_0)=0$ almost surely; on this event, fix $p\in A_0$ and $p_n\in A_{n,0}$. Every subsequential limit of $p_n$ lies in $A_0$, hence at $d$-distance zero from $p$, so we have $d(p_n,p)\to0$ by compactness of $T$ and continuity of $d$; by the triangle inequality, $d(t,p)$ does not depend on the choice of $p\in A_0$. For $t\in A_\eps$ maximizing $d(t,p)$ and $t_n\in T_n$ tending to $t$, we get $d(t,p)=\lim_nd(t_n,p_n)\le\liminf_n\max\{d(s,p'):s\in A_{n,\eps+\delta}, p'\in A_{n,0}\}$.

\emph{Step 6 (Increments).} Put $\Gamma_n=\max_{A_{n,\eps+\delta}}\widetilde G-\max_{A_{n,0}}\widetilde G$. For $t\in A_\eps$ maximizing $\widetilde G$ over $A_\eps$ and $t_n\in T_n$ tending to $t$, we have $\liminf_n\max_{A_{n,\eps+\delta}}\widetilde G\ge\lim_n\widetilde G(t_n)=\max_{A_\eps}\widetilde G$, while $\limsup_n\max_{A_{n,0}}\widetilde G\le\max_{A_0}\widetilde G$ because every subsequential limit of points of $A_{n,0}$ lies in $A_0$. Hence we have $0\le\Gamma_\eps\le\liminf_n\Gamma_n$. This proves the inequalities and the uniqueness assertion.
\end{proof}

\begin{remark}\label{rem:fromsharp}
Inequality~\eqref{eq:processsharp} alone does not give the uniqueness assertion, since $\Gamma_0\equiv0$; the proof obtains uniqueness from \eqref{eq:processbound} at $\eps=0$. Given uniqueness, however, \eqref{eq:processsharp} implies the first terms of the minima in \eqref{eq:processbound} and \eqref{eq:processradius}, by the argument of Lemma~\ref{lem:geometry} applied pathwise. Let $A_\theta(g)$ be as in the proof of Theorem~\ref{thm:process}, write $\E_{\widetilde G}$ for expectation over $\widetilde G$ with $g$ fixed, and fix $g\in C(T)$ with $\diam_d(A_0(g))=0$ and $p\in A_0(g)$. For $s,t\in A_0(g)$, the variance of $\widetilde G(s)-\widetilde G(t)$ is $d(s,t)^2=0$. Applying this to a countable dense subset of $A_0(g)$ and using the continuity of the paths, we see that $\widetilde G$ is almost surely constant on $A_0(g)$. Hence, almost surely, $\Gamma_\eps(g)=\max_{A_\eps(g)}\widetilde G-\max_{A_0(g)}\widetilde G$ equals $\max_{t\in A_\eps(g)}\{\widetilde G(t)-\widetilde G(p)\}$, and the same holds with $-\widetilde G$ in place of $\widetilde G$. The proof of Lemma~\ref{lem:geometry}, with $\ip Ux$ replaced by $\widetilde G(t)-\widetilde G(p)$ and with the symmetry of the law of $\widetilde G$ in place of that of $U$, gives
\[
 \diam_d(A_\eps(g))^2\le4\E_{\widetilde G}\Gamma_\eps(g)^2,\qquad
 \max_{t\in A_\eps(g)}d(t,p)^2\le2\E_{\widetilde G}\Gamma_\eps(g)^2.
\]
Since $d(t,p)$ does not depend on $p\in A_0(g)$, $\diam_d(A_0)=0$ almost surely, $X$ is independent of $\widetilde G$, and $(g,\tilde g)\mapsto\max_{A_\eps(g)}\tilde g-\max_{A_0(g)}\tilde g$ is Borel by the upper semicontinuity argument of Step~1 above, we may take $g=X$, for which $\Gamma_\eps(X)=\Gamma_\eps$, and integrate over $X$ by Tonelli's theorem. We obtain
\[
 \E\diam_d(A_\eps)^2\le4\E\Gamma_\eps^2,\qquad
 \E\max_{\substack{t\in A_\eps\\p\in A_0}}d(t,p)^2\le2\E\Gamma_\eps^2.
\]
\end{remark}

\section{Proof of the optimality of the constants}
\label{app:optimality}

\begin{proof}[Proof of Proposition~\ref{prop:sharpness}]
\emph{Step 1 (Diameter).} Take $f=0$, $\sigma=1$, $T=\{t\in\R^2:\norm t\le1\}$, and $G(t)=\ip Zt$ with $Z\sim N(0,I_2)$. Then $d$ is Euclidean, and we have $\E\sup_TG=\E\norm Z$ and $\norm Z>0$ almost surely. Write $e_1=(1,0)$, $e_2=(0,1)$, and rotate coordinates to obtain $Z=\norm Z e_1$. We have
\[
 A_\eps=\{x:\norm x\le1,\ x_1\ge1-\eta\},\qquad
 \eta=\eps/\norm Z.
\]
If $0<\eta<1$, we obtain
\[
 \norm{x-(1-\eta)e_1}^2=\norm x^2-2(1-\eta)x_1+(1-\eta)^2\le1-(1-\eta)^2=2\eta-\eta^2\qquad(x\in A_\eps).
\]
Hence we have $\diam(A_\eps)^2\le4(2\eta-\eta^2)$, with equality since $(1-\eta)e_1\pm\sqrt{2\eta-\eta^2} e_2\in A_\eps$. If $\eta\ge1$, then we have $\pm e_2\in A_\eps$ and $\diam(A_\eps)=2=\diam(T)$. This gives
\[
 \diam(A_\eps)^2=
 \begin{cases}
  8\eps/\norm Z-4\eps^2/\norm Z^2,&\eta<1,\\
  4,&\eta\ge1.
 \end{cases}
\]
For the Rayleigh density $s\mathrm e^{-s^2/2}$ of $\norm Z$ on $(0,\infty)$, we have $\E\norm Z^{-1}=\int_0^\infty\mathrm e^{-s^2/2} ds=\sqrt{\pi/2}$ and, integrating by parts, $\E\norm Z=\int_0^\infty s^2\mathrm e^{-s^2/2} ds=\sqrt{\pi/2}$.
Since $0\le\diam(A_\eps)^2/\eps\le8/\norm Z$ and the quotient tends almost surely to $8/\norm Z$, dominated convergence yields
\[
 \lim_{\eps\downarrow0}
 \frac{\E\diam(A_\eps)^2}{\eps \E\norm Z}
 =\frac{8\E\norm Z^{-1}}{\E\norm Z}=8.
\]

\emph{Step 2 (Distance from the maximizer).} Take $f=0$, $\sigma=1$, $T=[0,1]$, and $G(t)=Zt$ with scalar $Z\sim N(0,1)$; then we have $d(s,t)=\abs{s-t}$. Put $\ell=\min\{1,\eps/\abs Z\}$. The sets are $A_0=\{1\}$ and $A_\eps=[1-\ell,1]$ on $\{Z>0\}$, and $A_0=\{0\}$ and $A_\eps=[0,\ell]$ on $\{Z<0\}$. Hence we have
\[
 \max_{\substack{t\in A_\eps\\p\in A_0}}\abs{t-p}
 =\ell,\qquad
 \E\sup_TG=\E Z_+=\int_0^\infty z\ph(z) dz=\ph(0).
\]
As $\eps\downarrow0$, splitting at $\abs Z=\eps$ and substituting $z=\eps u$ gives
\begin{align*}
 \frac{\E\ell^2}{\eps}
 &=\frac2\eps\int_0^\eps\ph(z) dz
   +2\eps\int_\eps^\infty\frac{\ph(z)}{z^2} dz\\
 &=2\int_0^1\ph(\eps u) du
   +2\int_1^\infty\frac{\ph(\eps u)}{u^2} du
 \longrightarrow4\ph(0).
\end{align*}
Dominated convergence applies with bounds $\ph(0)$ on $[0,1]$ and $\ph(0)u^{-2}$ on $[1,\infty)$. This proves
\[
 \lim_{\eps\downarrow0}
 \frac{\E\max_{t\in A_\eps, p\in A_0}\abs{t-p}^2}
      {\eps \E\sup_TG}=4.
\]

\emph{Step 3 (Increments of an independent copy).} In the example of Step 2, let $\widetilde G(t)=Ut$ with scalar $U\sim N(0,1)$ independent of $Z$. By the form of $A_0$ and $A_\eps$ in Step 2, we have $\Gamma_\eps=\ell (-U)_+$ on $\{Z>0\}$ and $\Gamma_\eps=\ell U_+$ on $\{Z<0\}$. This gives $\E[\Gamma_\eps^2\mid Z]=\ell^2/2$ and $\E\Gamma_\eps^2=\frac12\E\ell^2$, and Step 2 yields $\E\Gamma_\eps^2/(\eps \E\sup_TG)\to2$ as $\eps\downarrow0$.

\emph{Step 4 (Fixed $\eps$).} Fix $\eps>0$ and $c>0$, and replace $T$ by $cT=\{ct:t\in T\}$ in the examples of Steps 1 and 2. Writing $A_\eps^T$ and $\Gamma_\eps^T$ to indicate the index set, we have $A_\eps^{cT}=cA_{\eps/c}^T$, $\Gamma_\eps^{cT}=c \Gamma_{\eps/c}^T$, and $\E\sup_{cT}G=c \E\sup_TG$, since $G$ and $\widetilde G$ are linear. Hence the three ratios of Steps 1--3 for $(cT,\eps)$ equal those for $(T,\eps/c)$, and they tend to $8$, $4$, and $2$ as $c\to\infty$. This proves the optimality assertions.
\end{proof}

\begin{proof}[Proof of \eqref{eq:balllimits}]
\emph{Step 1 (Diameter and distance from the maximizer).} Recall that $T$ is the Euclidean unit ball of $\R^k$ with $k\ge2$, $G(t)=\ip Zt$ with $Z\sim N(0,I_k)$, $f=0$, and $\sigma=1$; then $d$ is Euclidean and we have $\E\sup_TG=\E\norm Z$. Since the density of $\norm Z$ is proportional to $s^{k-1}\mathrm e^{-s^2/2}$, integration by parts gives $\E\norm Z^{-1}=\E\norm Z/(k-1)<\infty$. Let $w=Z/\norm Z$ be the maximizer and $\eta=\eps/\norm Z$. As in Step~1 of the preceding proof, we have $\diam(A_\eps)^2=8\eps/\norm Z-4\eps^2/\norm Z^2$ when $\eta\le1$. In this case we also have $\max_{t\in A_\eps}\norm{t-w}^2=2\eta$, since $\norm{t-w}^2=\norm t^2-2\ip tw+1\le2\eta$ on $A_\eps$, with equality when $\norm t=1$ and $\ip tw=1-\eta$. When $\eta>1$, both are at most $\diam(T)^2=4<4\eps/\norm Z$. Hence both are at most $8\eps/\norm Z$ for every $\eps>0$, and dominated convergence yields the first two limits in \eqref{eq:balllimits}.

\emph{Step 2 (Increments).} Let $\widetilde G(t)=\ip Ut$ with $U\sim N(0,I_k)$ independent of $Z$, and put $S=\ip Uw$ and $U_\perp=U-Sw$. If $\eta\le1$, every point of $A_\eps$ has the form $t=rw+y$ with $y\perp w$, $1-\eta\le r\le1$, and $\norm y^2\le1-r^2\le2\eta$. Hence $\Gamma_\eps=\max_{t\in A_\eps}\ip U{t-w}$ satisfies
\[
 \sqrt{2\eta-\eta^2} \norm{U_\perp}-\abs S\eta
 \le\Gamma_\eps
 \le\sqrt{2\eta} \norm{U_\perp}+\abs S\eta.
\]
The lower bound holds because $\ip U{t-w}=\sqrt{2\eta-\eta^2} \norm{U_\perp}-S\eta$ at the point of $A_\eps$ with $r=1-\eta$ and $y=\sqrt{2\eta-\eta^2} U_\perp/\norm{U_\perp}$. Thus we have $\Gamma_\eps^2/\eps\to2\norm{U_\perp}^2/\norm Z$ almost surely. On $\{\eta\le1\}$ we have $\Gamma_\eps^2/\eps\le(4\norm{U_\perp}^2+2S^2)/\norm Z$; on $\{\eta>1\}$ we have $\Gamma_\eps\le2\norm U$, which gives $\E[\eps^{-1}\Gamma_\eps^2\1_{\{\norm Z<\eps\}}]\le4k\eps^{-1}\Pbb\{\norm Z<\eps\}\to0$. Since $\E[\norm{U_\perp}^2\mid Z]=k-1$ and $\E\norm Z^{-1}<\infty$, dominated convergence gives $\lim_{\eps\downarrow0}\eps^{-1}\E\Gamma_\eps^2=2(k-1)\E\norm Z^{-1}=2\E\sup_TG$.
\end{proof}

\section{Proof of the drift-adaptive bounds}
\label{app:adaptive}

\begin{proof}[Proof of Proposition~\ref{prop:adaptive}]
\emph{Step 1 (Finite index sets).} In the setting of Proposition~\ref{thm:finite} and the notation of its proof, fix $\lambda>0$ and $s>0$. For $x\in C_{N,s}(z)$, the definition of $F_N^*$ gives $F_N(z')\ge\ip x{z'}-F_N^*(x)\ge F_N(z)+\ip x{z'-z}-s$ for all $z'\in\R^k$. With $z'=(1+\lambda)z$ and \eqref{eq:reg-conjugate}, we obtain
\begin{align*}
 \ip zx&\le\frac{F_N((1+\lambda)z)-F_N(z)+s}{\lambda}\\
 &\le\frac{F((1+\lambda)z)-F(z)+s}{\lambda}
   +\frac{\log m}{\lambda N}+\frac{(1+\lambda)^2-1}{2\lambda N}\norm z^2.
\end{align*}
The support points $q_s(Z,U)$ in \eqref{eq:supportidentity}, for $H=F_N^*$, lie in $C_{N,s}(Z)$. With this estimate in place of \eqref{eq:widthbudget}, the argument leading to \eqref{eq:regbounds}, that is, integration over $0<s\le\theta$ and passage to the limit superior as $N\to\infty$, gives
\begin{align*}
 \limsup_N\E h_{N,\theta}(Z,U)^2
 &\le\frac{2\theta \E[F((1+\lambda)Z)-F(Z)]+\theta^2}{\lambda},\\
 \limsup_N\E\diam(C_{N,\theta}(Z))^2
 &\le\frac{8\theta \E[F((1+\lambda)Z)-F(Z)]+4\theta^2}{\lambda},\\
 \limsup_N\E\max_{x\in C_{N,\theta}(Z)}\norm{x-\nabla F_N(Z)}^2
 &\le\frac{4\theta \E[F((1+\lambda)Z)-F(Z)]+2\theta^2}{\lambda}.
\end{align*}
The remainder of the proof of Proposition~\ref{thm:finite} is unchanged. Hence its three conclusions hold with $8\eps \E\max_i\ip Z{v_i}$, $4\eps \E\max_i\ip Z{v_i}$, and $2\eps \E\max_i\ip Z{v_i}$ replaced by the corresponding right-hand sides above, evaluated at $\theta=\eps$.

\emph{Step 2 (Compact index sets).} For finite $T_n\subset T$, the reduction in Appendix~\ref{sec:process} with $a_i=f(t_i)/\sigma$ and $\eps/\sigma$ in place of $\eps$ gives $\sigma \E F(cZ)=\E\max_{T_n}(f+c\sigma G)$ for $c\in\{1,1+\lambda\}$. Hence \eqref{eq:adaptivediam}--\eqref{eq:adaptivesharp} hold for $T_n$ with $\Delta_n(\lambda)=\E\max_{T_n}(f+(1+\lambda)\sigma G)-\E\max_{T_n}(f+\sigma G)$ in place of $\Delta(\lambda)$. For the nets $T_n$ constructed in Appendix~\ref{sec:process}, continuity of the paths gives $\max_{T_n}(f+c\sigma G)\uparrow\max_T(f+c\sigma G)$ for $c\in\{1,1+\lambda\}$. These maxima are dominated by $\norm f_\infty+(1+\lambda)\sigma\norm G_\infty$, so we have $\Delta_n(\lambda)\to\Delta(\lambda)$. The Fatou arguments of Appendix~\ref{sec:process} then give \eqref{eq:adaptivediam}--\eqref{eq:adaptivesharp} for $T$.

\emph{Step 3 (Limits of $\Delta(\lambda)/\lambda$).} Fix a continuous path and put $\Phi(\lambda)=\max_T(X+\lambda\sigma G)$ for $\lambda\in\R$; then we have $\Delta(\lambda)=\E[\Phi(\lambda)-\Phi(0)]$ for $\lambda>0$. As a maximum of affine functions of $\lambda$, $\Phi$ is convex; hence $(\Phi(\lambda)-\Phi(0))/\lambda$ is nondecreasing in $\lambda>0$, and so is $\Delta(\lambda)/\lambda$. Since $T$ is compact and $(t,\lambda)\mapsto X(t)+\lambda\sigma G(t)$ is continuous and affine in $\lambda$, Danskin's theorem~\cite{Danskin1966} gives the right derivative $\sigma\max_{t\in A_0}G(t)$ of $\Phi$ at $0$; therefore we have $(\Phi(\lambda)-\Phi(0))/\lambda\downarrow\sigma\max_{t\in A_0}G(t)$ as $\lambda\downarrow0$. Since $\abs{\Phi(\lambda)-\Phi(0)}\le\lambda\sigma\norm G_\infty$, dominated convergence gives $\Delta(\lambda)/\lambda\to\sigma\E\max_{A_0}G$. As $\lambda\to\infty$, we have $\Phi(\lambda)/\lambda=\max_T(X/\lambda+\sigma G)\to\sigma\sup_TG$ for every path. Moreover, we have $\abs{\Phi(\lambda)}/\lambda\le\norm f_\infty+2\sigma\norm G_\infty$ for $\lambda\ge1$ and $\Phi(0)/\lambda\to0$. Dominated convergence therefore gives $\Delta(\lambda)/\lambda\to\sigma\E\sup_TG$.

\emph{Step 4 (Small $\eps$).} Given $\eps'>0$, choose $\lambda>0$ with $\Delta(\lambda)/\lambda\le\sigma\E\max_{A_0}G+\eps'$. By \eqref{eq:adaptivediam}, we have $\eps^{-1}\E\diam_d(A_\eps)^2\le8\Delta(\lambda)/(\lambda\sigma^2)+4\eps/(\lambda\sigma^2)$, whose limit superior as $\eps\downarrow0$ is at most $8(\E\max_{A_0}G+\eps'/\sigma)/\sigma$. Letting $\eps'\downarrow0$ proves the first inequality in \eqref{eq:adaptivelimit}; the other two follow from \eqref{eq:adaptiveradius} and \eqref{eq:adaptivesharp} in the same way.
\end{proof}

\bibliographystyle{alpha}

\bibliography{references}

\end{document}